\documentclass[12pt]{amsart}

\usepackage[margin=1in,marginparwidth=0.8in, marginparsep=0.1in]{geometry}

\pdfoutput=1

\usepackage{soul}

\usepackage{amsfonts,amssymb,latexsym,amsmath,graphicx}

\usepackage{hyperref}

\usepackage[all, knot, poly]{xy}
\xyoption{knot}

\usepackage{times}
\usepackage{mathrsfs}

\newtheorem{lemma}{Lemma}

\newtheorem{theorem}[lemma]{Theorem}

\theoremstyle{definition}

\theoremstyle{remark} 
\newtheorem*{remark}{Remark}

\newtheorem{question}[lemma]{Question}

\newcommand{\C}{\mathbb{C}}

\newcommand{\cA}{\mathcal{A}}

\newcommand{\cH}{\mathcal{H}}

\newcommand{\cL}{\mathcal{L}}

\newcommand{\bX}{\mathbf{X}}

\newcommand{\Ho}{\mathrm{H}}
\newcommand{\SH}{\mathrm{SH}}

\DeclareMathOperator{\Hom}{Hom}

\DeclareMathOperator{\colim}{colim}
\DeclareMathOperator{\Fuk}{Fuk}

\usepackage{tikz-cd}

\author{Vivek Shende}

\title{An algebraic approach to the algebraic Weinstein conjecture}

\begin{document}

\begin{abstract}
How does one measure the failure of Hochschild homology to commute with colimits? 

Here I relate this question to 
a major open problem about dynamics in contact manifolds --- 
the assertion that Reeb orbits exist and are detected by symplectic homology. 
More precisely, I show that for polarizably Weinstein fillable contact manifolds,
said property is equivalent to the failure of Hochschild homology to commute with
certain colimits of representation categories of tree quivers. 

So as to be intelligible to algebraists, I try to include or black-box as much of the geometric
background as possible.  
\end{abstract}

\maketitle

\thispagestyle{empty}

Existence of closed geodesics on a compact Riemannian manifold $M$ is guaranteed for topological reasons \cite{LF, 
Klingenberg}.  
Let us recall the argument. 
Morse theory tells us that the homology of the free loop space $\cL M = \mathrm{Maps}(S^1 \to  M)$ 
can be computed by a complex generated by geodesics.  The trivial loops 
contribute a subcomplex computing the homology of the original manifold, so 
there must be nontrivial loops unless $\Ho_\bullet(M) \to \Ho_\bullet(\cL M)$ is an isomorphism. 
It is obviously not an isomorphism unless $M$ is simply connected, and in this case
we study the based loop space $\Omega M := \mathrm{Maps}((S^1, 0) \to (M, m))$
and the fibration $\Omega M \to \cL M \to M$.  As this is split by constant loops
$M \to \cL M$, we find $\pi_k(\cL M) = \pi_k(M) \oplus \pi_{k+1}(M)$, so by Hurewicz
the first nontrivial homotopy group $\pi_{k+1}(M)$ contributes nontrivially to $\Ho_k(\cL M)$, while 
$\Ho_k(M)$ vanishes.  In fact, in the simply connected case, one can obtain rather 
more refined information \cite{VP-Sullivan}. 

We would like to think of the map $\Ho_\bullet(M) \to \Ho_\bullet(\cL M)$ as arising from 
the following  local-to-global construction. 
On $M$, we can
consider the constant  cosheaf of spaces, which we denote $\Omega$.  By definition, 
this assigns a point to any contractible
open set, and sends covers to colimits.\footnote{We abuse terminology to respect intuition and write `homology' 
when we mean that we think of the complex as an object in the derived category.  
Similarly we write `$=$' to mean e.g. identified by a canonical
quasi-isomorphism, etc.  We always work with $\infty$-categories, etc., and (pre)sheaves or cosheaves should be understood
accordingly.}  
By the ($\infty$) van Kampen theorem, 
$\Omega(U)$ is the path groupoid of $U$, explaining the notation.  Composition
of loops gives $\Ho_\bullet(\Omega(M))$ the structure of a ring, and its 
Hochschild homology gives the homology of the free loop space \cite{Burghelea-Fiedorowicz, Goodwillie}: 
$\Ho \Ho_\bullet (\Ho_\bullet(\Omega M)) = \Ho_\bullet(\cL M)$.  The inclusion 
of constant loops is:  
\begin{equation} \label{failure} \Ho_\bullet(M) = \colim_U \Ho \Ho_\bullet (\Ho_\bullet(\Omega (U)))
\to \Ho \Ho_\bullet (\colim_U \Ho_\bullet(\Omega (U))) = \Ho \Ho_\bullet (\cL M) \end{equation}
The following is the first example of the problem we are interested in: 

\begin{question}
Can the failure of (\ref{failure}) to be an isomorphism be seen in terms of some general machinery 
measuring the failure of Hochschild homology to commute with homotopy colimit?
\end{question}

\vspace{2mm} 
We turn to contact geometry. 
The contact-geometric formulation of the geodesic flow is the following.  On a cotangent bundle $T^*M$,
there is the tautological 1-form $\lambda$, which at a given covector $\xi$ is the function on tangent vectors
given by $\xi$.  Fixing a metric on $M$, we may restrict $\lambda$ to the cosphere bundle $S^*M$; here it is
{\em contact}, meaning that $\lambda \wedge (d\lambda^{n-1})$ is nowhere vanishing.  The {\em Reeb}
vector field $R$ on $S^*M$ is characterized by lying in the kernel of $d \lambda$ and normalized by $\lambda(R) = 1$.  
Its flow is naturally identified with the geodesic flow.  

More generally, the same formulas define the Reeb flow for any contact form on any (odd dimensional) manifold $V$.  
The Weinstein conjecture asserts the existence of a closed trajectory \cite{Weinstein-Conjecture}.  It has long been
a central problem in contact geometry 
\cite{Viterbo-WeinsteinR2n, 
Hofer-Viterbo-Cotangent, 
Hofer-Viterbo-Spheres, 
Hofer-Weinstein-3}, 
and known in general only in dimension 3, by an argument whose ingredients have no known analogue in higher dimension
\cite{Taubes-W}.  The result is also known for flexible contact structures in general \cite{Albers-Hofer}; we will
be interested here in what is in some sense the opposite setting \cite{Niederkruger}, of Weinstein fillable contact manifolds.  

It is natural to try and generalize the Morse theoretic approach to geodesics to the study of Reeb orbits.  That is,
one wants a complex generated by orbits, so that nonvanishing of the homology 
groups implies the existence of orbits.  The reason to impose a differential is to provide invariance
under deformations: these homologies depend only on $\mathrm{ker}\, \lambda$ rather
than $\lambda$ itself, i.e. on the {\em contact structure} rather than the {\em contact form}.\footnote{Taubes's celebrated 
work in 3 dimensions \cite{Taubes-W} also can be understood in terms of 
such a homology theory \cite{Taubes-ECH}, though the theory he constructs (embedded contact homology) 
turns out to in fact be completely independent of $\lambda$, and in fact isomorphic to a topological invariant
(Seiberg-Witten homology) known to be sufficiently nontrivial to imply existence of orbits.  
Unfortunately, in higher dimensions, there are no known analogues of the embedded contact homology or the 
Seiberg-Witten homology.}

Just as the cosphere bundle bounds the codisk bundle, we may ask that some general contact 
$(V, \lambda)$ is the boundary of some $W$ to which $\lambda$ extends, and over which $(d \lambda)^n$ 
is everywhere nonvanishing.  We also ask that the `Liouville' vector field  $Z$ characterized by 
$d\lambda(Z, \cdot) = \lambda$ points out at the boundary; in the case of the cotangent bundle this vector 
field is radial in the fiber directions.  Such $W$ are called Liouville domains, and determine a 
symplectic cohomology $\SH^\bullet(W)$, which may be taken to be 
generated by the Reeb orbits of $(V, \lambda)$ and the critical points of a Morse function on $W$ \cite{Viterbo-Functors1}.  
(See \cite{Seidel} for a leisurely introduction to Liouville manifolds and symplectic cohomology.)   The differential 
is such that these Morse critical points form a subcomplex on which the differential is the Morse differential, giving an
exact triangle.\footnote{There are many differing convention for the grading of symplectic cohomology,
and also for which item to call symplectic cohomology and which symplectic homology.  We follow \cite{Abouzaid-Generation, GPS1}.} 
\begin{equation} \label{viterbo triangle} \Ho^\bullet(W) \to \SH^\bullet(W) \to \widetilde{\SH}{}^\bullet(W) \xrightarrow{[1]} \end{equation}
In particular, if $\Ho^\bullet(W) \to \SH^\bullet(W)$ fails to be an isomorphism, then the Weinstein conjecture holds for $V$, 
as  $\widetilde{\SH}{}^\bullet(W)$, which is generated by the Reeb orbits of $V$, must be nonzero.  

Viterbo's {\em algebraic Weinstein conjecture} is the assertion that $\widetilde{\SH}{}^\bullet(W)$ always
detects an orbit.  It is not typically easy to compute $\SH^\bullet(W)$.   But when $W = T^*M$, one knows: 
\begin{theorem} \cite{Viterbo-ICM, Viterbo-Functors2}
There is an isomorphism $\SH^{\bullet + n} (T^*M) \cong \Ho_{- \bullet}(\cL M)$.\footnote{When $M$ is not spin, 
it is necessary to twist one side or the other by a local system \cite{Kragh}.} 
\end{theorem}
\noindent This isomorphism has seen many further developments; see e.g. \cite{Ciliebak-Latschev, Abouzaid-Viterbo}. 

\vspace{2mm}

The composition
$\Ho^{\bullet + n} (T^*M) \to \SH^{\bullet + n} (T^*M) \cong \Ho_{- \bullet}(\cL M)$ is identified
with the inclusion of constant loops $\Ho_{-\bullet}(M) \to \Ho_{-\bullet}(\cL M)$ under Poincar\'e duality \cite[Lem. 3.6]{Abouzaid-Cotangent}. 
So, the symplectic homology detects geodesics in essentially the same way as the Morse homology of the loop space
did.  However already in this case it does more: it shows that a contact level of $T^*M$, not necessarily 
the unit cosphere bundle for any Riemannian metric, will also have Reeb orbits.

\vspace{2mm}
A class of Liouville domains including but rather more general than codisk bundles are the {\em Weinstein} 
domains.\footnote{Weinstein manifolds and the Weinstein conjecture have the same eponym, but a priori no other relation.} 
By definition, these are those for which the Liouville vector field is gradient-like for a Morse-Bott 
function. In this case the critical points of $Z$ have index $\le \dim W /2$, and union
of descending level sets is a singular isotropic subset termed the {\em core} or {\em skeleton}.  Stein domains from
complex analysis are Weinstein when viewed as symplectic manifolds, and conversely any Weinstein
domain is deformation equivalent to a Stein domain \cite{Cieliebak-Eliashberg}. 

Weinstein domains for which the indices of critical points are $< \dim W /2$ are said to be {\em subcritical}, and
for these it is known that $\SH^\bullet(W) = 0$; in particular, the Weinstein conjecture holds for their contact boundaries
\cite{Viterbo-Functors2}.  Simple examples: the ball; the cotangent bundle of a noncompact manifold.
Beyond these, the Weinstein conjecture
is not known for contact boundaries of Weinstein domains in any reasonable generality, 
and it would be a major advance to establish it. 

\vspace{2mm}

One available tool for computing symplectic homology is the open-closed
morphism from the Hochschild homology of the wrapped Fukaya category\footnote{\cite{Seidel-book, GPS1} contain the relevant definitions.  
We will soon cite some results which compute the Fukaya category in all relevant cases, so 
 the gist of the article will not be lost to the reader with no idea what the Fukaya category is.} 
\cite{FOOO, Seidel-Deformations, Seidel-Hochschild, BEE, Abouzaid-Generation, GPS1}:
\begin{equation} \Ho \Ho_\bullet (\Fuk(W)) \to \SH^{\bullet+n}(W) \end{equation} 
By either \cite{Abouzaid-Generation, Ganatra} or \cite{BEE, Ekholm-Lekili, Ekholm-Surgery} plus a `generation' result 
\cite{CDGG, GPS2}, this morphism is by now known to be an isomorphism for Weinstein domains.  

For cotangent bundles, there is an object (the cotangent fiber) of $F \in \Fuk(T^* M)$, which 
generates the category and for which $\Hom(F, F) = \Ho_{-\bullet}(\Omega M)$.  Thus the open-closed
map induces 

$$ \Ho \Ho_\bullet  ( \Ho_{-\bullet}(\Omega M) ) =  \Ho \Ho_\bullet ( \Hom(F, F)) \to \SH^{\bullet + n} (T^*M) = \Ho_{- \bullet}(\cL M) $$
This is the same as the corresponding such morphism mentioned above.\footnote{I am not 
certain whether this follows from \cite{Abouzaid-Cotangent}, but in any case it certainly does from \cite{GPS1, GPS2, GPS3}.} 

How can the open closed map help us?  At first, it does not look promising.  
We have not said what the Fukaya category is, but its definition involves
the same sort of geometrical structures as are involved in symplectic homology.  
On top of this we have now added the nontrivial step of taking Hochschild homology.  
However, just as $\Omega M$ has better local-to-global behavior than $\cL M$, we also have
(as anticipated by \cite{Kontsevich}): 

\begin{theorem} \cite{GPS1, GPS2, GPS3, Shende-Microlocal, Nadler-Shende}
The Fukaya category of a Weinstein manifold is the global sections of a constructible cosheaf of categories over
the skeleton.  Moreover, this cosheaf is isomorphic to the cosheaf of microlocal sheaves.\footnote{Strictly 
speaking, what is presently in the literature requires that the Weinstein manifold is `stably polarizable'.  
It is known to experts how to remove this hypothesis; on the other hand  we will impose it later for different reasons.} 
\end{theorem}

We have not said what microlocal sheaves are and it will not be relevant; 
but for a definition see \cite{Nadler-Shende}, which is
built on the technology of \cite{Kashiwara-Schapira}. 
What is relevant is that microlocal sheaves are in principle combinatorial-topological in nature, 
but in practice the stalks of the above cosheaf
may be complicated categories at complicated singularities of the skeleton.  When the skeleton
is smooth, the cosheaf is locally constant with stalk the category of chain complexes.  

For $W = T^*M$,
the cosheaf is simply the path groupoid $\Omega$ (twisted by a local system if $M$ is not spin). 
More generally, Nadler found an explicit collection of so-called `arboreal' singularities \cite{Nadler-Arboreal-1, Nadler-Arboreal-2}, 
with the property that: 

\begin{theorem} \cite{Nadler-Arboreal-1}
When the skeleton is arboreal, the cosheaf of microlocal sheaves has stalks given by representation
categories of tree quivers.  The cogenerization morphisms are explicit.  
\end{theorem} 

\begin{figure}
\begin{center}
\includegraphics[scale=0.35]{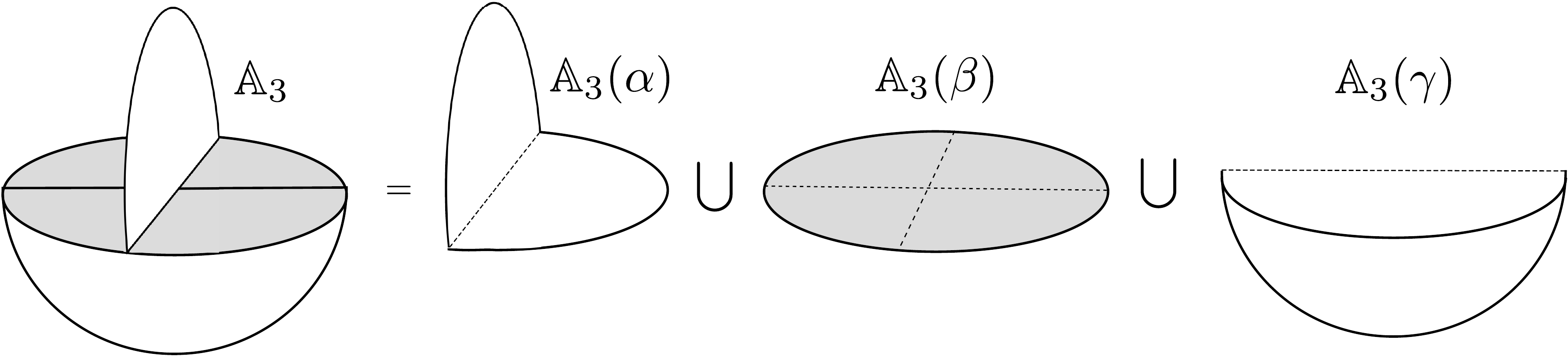}
\end{center}
\caption{  The arboreal A3 singularity, and how to glue it from smooth spaces.  All arboreal
singularities admit analogous gluing descriptions. }
\label{A3}
\end{figure}

When $\dim W = 2$, the skeleton of $W$ will be a (ribbon) graph, and `arboreal' essentially 
amounts to asking that the graph is trivalent.  The cosheaf $\cA$ will assign the category of 
chain complexes at smooth points, and the category of exact triangles (aka
$\mathrm{Perf}(\bullet \to \bullet)$) at the trivalent vertices, with the obvious cogenerization morphisms. 
This case was studied in \cite{Dyckerhoff-Kapranov}. 
For $\dim W = 4$, the skeleton is two dimensional;  the topology of the 
typical new kind of singularity is depicted in Figure \ref{A3}. 
More geometric pictures can be found in \cite{Starkston-Arboreal, AgEN}. 

By an arboreal space, we will mean a pair $(\bX, \cA)$ of a space and a constructible cosheaf, which are locally given
by Nadler's explicit models.  In particular, the stalks are representation categories of tree quivers.  
An explicitly combinatorial exposition of this notion can be found in \cite{Shende-Takeda}.  

\vspace{2mm}

Some further
restrictions on the local models lead to the notion of {\em positive arboreal space} \cite{AgEN}.   
These provide skeleta for a large class of Weinstein manifolds: 

\begin{theorem} \cite{AgEN}
There is an equivalence of categories between: 
\begin{itemize}
\item Stably polarized Weinstein manifolds \& their homotopies
\item Positive arboreal spaces \& their concordances
\end{itemize}
\end{theorem}

Here, a stable polarization is a choice of Lagrangian sub-bundle of $TW \oplus \C^n$ for some $n$.   

\begin{remark} 
It is
not clearly understood what weaker hypothesis corresponds to taking all arboreal spaces; in dimension
4, none is needed \cite{Starkston-Arboreal}, but this is apparently not the case in general \cite{AgEN}.  It is
expected that this can be repaired by adding some further explicit list of presently unknown singularities. 
There are also other tricks for reducing problems to the stably polarized case, like \cite[Sec. 10]{Nadler-Shende}. 
In any case, I expect that any technique which works for stably polarized Weinstein manifolds should work for
Weinstein manifolds in general. 
\end{remark} 

In some sense we have already arrived at a reduction to algebra: one could try and develop tools for computing
the colimit of categories giving $\cA(\bX)$ or its Hochschild homology
$\Ho\Ho_\bullet(\cA(\bX))$.  Indeed, as far
as anyone knows, $\SH^\bullet(W)$ is always either infinite dimensional, or zero.  If one could show this
zero or infinite property for $\Ho\Ho_\bullet(\cA(\bX))$, the algebraic Weinstein conjecture would follow.  Or if one could show that the 
natural circle action on any nonvanishing $\Ho\Ho_\bullet(\cA(\bX))$ is nontrivial, the result would again follow.  

Here we want to point out that in fact it is possible to directly generalize Equation (\ref{failure}) and make 
direct contact with Equation (\ref{viterbo triangle}).   Consider
a cosheaf of categories $\cA$ over a space $\bX$.  The Hochschild homologies form a precosheaf
$\Ho \Ho_\bullet(\cA)$ which is {\em not} generally a cosheaf, since Hochschild homology does not commute with colimit. 
We may cosheafify it and obtain a cosheaf 
$\cH \cH_\bullet(\cA)$.   There is a natural map 
\begin{equation}\Gamma(\bX, \cH \cH_\bullet(\cA)) \to \Ho \Ho_\bullet (\cA(\bX))\end{equation}
Because $\cA$ is constructible, the LHS is the colimit of the Hochschild homologies of the stalks of $\cA$, and
the RHS is the Hochschild homology of the colimit of the stalks.  

Moreover, $\cH \cH_\bullet(\cA)$ can be computed explicitly.  Indeed, for a tree quiver $T$, 
it is the case that $\Ho \Ho_\bullet (\mathrm{Perf}(T))$ is concentrated in degree zero,
and is a free module whose rank is the number of vertices of $T$.   This gives the 
stalks of $\cH \cH_\bullet(\cA)$, and in fact one can show: 

\begin{theorem} \label{local hochschild homology} \cite{Shende-Takeda}
When $(\bX, \cA)$ arises from the skeleton of a stably polarizable Weinstein manifold, 
$\cH \cH_{\bullet - n}(\cA)$ is the cosheaf of compactly supported cohomologies.  As $\bX$ is compact, 
we have $\Gamma(\bX, \cH \cH_{\bullet-n}(\cA)) \cong \Ho^\bullet (\bX)$. 
\end{theorem}

In fact, it is possible to show that the resulting map
\begin{equation}
\Ho^\bullet (\bX) \cong \Gamma(\bX, \cH \cH_{\bullet-n}(\cA))  \to \Ho \Ho_{\bullet - n} (\cA(\bX)) \to \SH^\bullet(W)
\end{equation}
agrees with the original $\Ho^\bullet (\bX) = \Ho^\bullet(W) \to \SH^\bullet(W)$.  This follows from 
 \cite[Eq. 1.7]{GPS1}
given that the local arboreal models are {\em Weinstein sectors}, as was shown in 
\cite{Shende-Arboreal}.  In some more detail: in \cite{Shende-Arboreal} 
it is shown that the nondegenerate arboreal sectors are {\em stopped}; hence the top row of \cite[Eq. 1.7]{GPS1}
consists of isomorphisms.  Degenerate arboreal singularities are obtained by stop removal;  the desired
commutativity descends using
the stop removal localization sequence and the fact that Hochschild homology sends localizations to exact
triangles. 

\vspace{2mm}
Putting all this together we have: 

\begin{theorem}
The algebraic Weinstein conjecture for contact manifolds with stably polarizable Weinstein filling is equivalent 
to the assertion that  $\Gamma(\bX, \cH \cH_\bullet(\cA)) \to \Ho \Ho_\bullet (\cA(\bX))$ is never an isomorphism for positive arboreal
spaces $(\bX, \cA)$. 
\end{theorem} 

\vspace{2mm}

We are left with the following: 

\begin{question}
What measures the failure of Hochschild homology to commute with colimits?
\end{question}

\vspace{2mm}
{\bf Acknowledgements.}  Beyond the evident intellectual debt that this work owes to the ideas of Claude Viterbo, 
it happens that I was inspired to write it by the upcoming celebration of his mathematics on the occasion of 
his sixtieth birthday.

\bibliographystyle{plain}
\bibliography{hochschild}

\begin{thebibliography}{10}

\bibitem{Abouzaid-Generation}
Mohammed Abouzaid.
\newblock A geometric criterion for generating the {F}ukaya category.
\newblock {\em Publications Math{\'e}matiques de l'IH{\'E}S}, 112:191--240,
  2010.

\bibitem{Abouzaid-Cotangent}
Mohammed Abouzaid.
\newblock A cotangent fibre generates the {F}ukaya category.
\newblock {\em Advances in Mathematics}, 228(2):894--939, 2011.

\bibitem{Abouzaid-Viterbo}
Mohammed Abouzaid.
\newblock Symplectic cohomology and {V}iterbo's theorem.
\newblock {\em arXiv:1312.3354}, 2013.

\bibitem{Albers-Hofer}
Peter Albers and Helmut Hofer.
\newblock On the {W}einstein conjecture in higher dimensions.
\newblock {\em arXiv:0705.3953}.

\bibitem{AgEN}
Daniel Alvarez-Gavela, Yakov Eliashberg, and David Nadler.
\newblock Positive arborealization of polarized {W}einstein manifolds.
\newblock {\em arXiv:2011.08962}.

\bibitem{BEE}
Fr{\'e}d{\'e}ric Bourgeois, Tobias Ekholm, and Yakov Eliashberg.
\newblock Effect of {L}egendrian surgery.
\newblock {\em Geometry \& Topology}, 16(1):301--389, 2012.

\bibitem{Burghelea-Fiedorowicz}
Dan Burghelea and Zbigniew Fiedorowicz.
\newblock Cyclic homology and algebraic k-theory of spaces—ii.
\newblock {\em Topology}, 25(3):303--317, 1986.

\bibitem{CDGG}
Baptiste Chantraine, Georgios~Dimitroglou Rizell, Paolo Ghiggini, and Roman
  Golovko.
\newblock Geometric generation of the wrapped {F}ukaya category of {W}einstein
  manifolds and sectors.
\newblock {\em arXiv:1712.09126}.

\bibitem{Cieliebak-Eliashberg}
Kai Cieliebak and Yakov Eliashberg.
\newblock {\em From {S}tein to {W}einstein and back: symplectic geometry of
  affine complex manifolds}.
\newblock American Mathematical Soc., 2012.

\bibitem{Ciliebak-Latschev}
Kai Cieliebak and Janko Latschev.
\newblock The role of string topology in symplectic field theory.
\newblock {\em New perspectives and challenges in symplectic field theory},
  49:113--146, 2009.

\bibitem{Dyckerhoff-Kapranov}
Tobias Dyckerhoff and Mikhail Kapranov.
\newblock Triangulated surfaces in triangulated categories.
\newblock {\em arXiv:1306.2545}.

\bibitem{Ekholm-Surgery}
Tobias Ekholm.
\newblock Holomorphic curves for {L}egendrian surgery.
\newblock {\em arXiv:1906.07228}.

\bibitem{Ekholm-Lekili}
Tobias Ekholm and Yanki Lekili.
\newblock Duality between {L}agrangian and {L}egendrian invariants.
\newblock {\em arXiv:1701.01284}.

\bibitem{LF}
Abram Fet and Lazar Lyusternik.
\newblock Variational problems on closed manifolds.
\newblock {\em Dokl. Akad. Nauk. SSSR}, 81:17--18, 1951.

\bibitem{FOOO}
Kenji Fukaya, Yong-Geun Oh, Hiroshi Ohta, and Kaoru Ono.
\newblock {\em Lagrangian intersection {F}loer theory: anomaly and
  obstruction}.
\newblock American Mathematical Soc., 2010.

\bibitem{Ganatra}
Sheel Ganatra.
\newblock Symplectic cohomology and duality for the wrapped {F}ukaya category.
\newblock {\em arXiv:1304.7312}.

\bibitem{GPS3}
Sheel Ganatra, John Pardon, and Vivek Shende.
\newblock Microlocal {M}orse theory of wrapped {F}ukaya categories.
\newblock {\em arXiv:1809.08807}.

\bibitem{GPS2}
Sheel Ganatra, John Pardon, and Vivek Shende.
\newblock Sectorial descent for wrapped {F}ukaya categories.
\newblock {\em arXiv:1809.03427}.

\bibitem{GPS1}
Sheel Ganatra, John Pardon, and Vivek Shende.
\newblock Covariantly functorial wrapped {F}loer theory on {L}iouville sectors.
\newblock {\em Publications math{\'e}matiques de l'IH{\'E}S}, 131(1):73--200,
  2020.

\bibitem{Goodwillie}
Thomas~G Goodwillie.
\newblock Cyclic homology, derivations, and the free loopspace.
\newblock {\em Topology}, 24(2):187--215, 1985.

\bibitem{Hofer-Weinstein-3}
Helmut Hofer.
\newblock Pseudoholomorphic curves in symplectizations with applications to the
  {W}einstein conjecture in dimension three.
\newblock {\em Inventiones mathematicae}, 114(1):515--563, 1993.

\bibitem{Hofer-Viterbo-Cotangent}
Helmut Hofer and Claude Viterbo.
\newblock The {W}einstein conjecture in cotangent bundles and related results.
\newblock {\em Annali della Scuola Normale Superiore di Pisa-Classe di
  Scienze}, 15(3):411--445, 1988.

\bibitem{Hofer-Viterbo-Spheres}
Helmut Hofer and Claude Viterbo.
\newblock The weinstein conjecture in the presence of holomorphic spheres.
\newblock {\em Communications on pure and applied mathematics}, 45(5):583--622,
  1992.

\bibitem{Kashiwara-Schapira}
Masaki Kashiwara and Pierre Schapira.
\newblock {\em Sheaves on Manifolds}.
\newblock Springer Science \& Business Media, 2013.

\bibitem{Klingenberg}
Wilhelm Klingenberg.
\newblock {\em Lectures on closed geodesics}.
\newblock Springer Science \& Business Media, 2012.

\bibitem{Kontsevich}
Maxim Kontsevich.
\newblock Symplectic geometry of homological algebra.
\newblock {\em available at the author’s webpage}, 2009.

\bibitem{Kragh}
Thomas Kragh.
\newblock The {V}iterbo transfer as a map of spectra.
\newblock {\em arXiv:0712.2533}.

\bibitem{Nadler-Arboreal-2}
David Nadler.
\newblock Non-characteristic expansions of {L}egendrian singularities.
\newblock {\em arXiv:1507.01513}.

\bibitem{Nadler-Arboreal-1}
David Nadler.
\newblock Arboreal singularities.
\newblock {\em Geometry \& Topology}, 21(2):1231--1274, 2017.

\bibitem{Nadler-Shende}
David Nadler and Vivek Shende.
\newblock Sheaf quantization in {W}einstein symplectic manifolds.
\newblock {\em arXiv:2007.10154}.

\bibitem{Niederkruger}
Klaus Niederkr{\"u}ger.
\newblock The plastikstufe -- a generalization of the overtwisted disk to
  higher dimensions.
\newblock {\em Algebraic \& Geometric Topology}, 6(5):2473--2508, 2006.

\bibitem{Seidel-Deformations}
Paul Seidel.
\newblock Fukaya categories and deformations.
\newblock {\em math/0206155}.

\bibitem{Seidel}
Paul Seidel.
\newblock A biased view of symplectic cohomology.
\newblock {\em Current developments in mathematics}, 2006(1):211--254, 2006.

\bibitem{Seidel-book}
Paul Seidel.
\newblock {\em Fukaya categories and Picard-Lefschetz theory}, volume~10.
\newblock European Mathematical Society, 2008.

\bibitem{Seidel-Hochschild}
Paul Seidel.
\newblock Symplectic homology as {H}ochschild homology.
\newblock {\em Algebraic geometry—Seattle 2005}, pages 415--434, 2009.

\bibitem{Shende-Arboreal}
Vivek Shende.
\newblock Arboreal singularities from {L}efschetz fibrations.
\newblock {\em arXiv:1809.10359}.

\bibitem{Shende-Microlocal}
Vivek Shende.
\newblock Microlocal category for {W}einstein manifolds via h-principle.
\newblock {\em arXiv:1707.07663}.

\bibitem{Shende-Takeda}
Vivek Shende and Alex Takeda.
\newblock Calabi-{Y}au structures on topological {F}ukaya categories.
\newblock {\em arXiv:1605.02721}.

\bibitem{Starkston-Arboreal}
Laura Starkston.
\newblock Arboreal singularities in {W}einstein skeleta.
\newblock {\em Selecta Mathematica}, 24(5):4105--4140, 2018.

\bibitem{Taubes-W}
Clifford~Henry Taubes.
\newblock The {S}eiberg--{W}itten equations and the {W}einstein conjecture.
\newblock {\em Geometry \& Topology}, 11(4):2117--2202, 2007.

\bibitem{Taubes-ECH}
Clifford~Henry Taubes.
\newblock Embedded contact homology and {S}eiberg--{W}itten {F}loer cohomology
  {I}.
\newblock {\em Geometry \& Topology}, 14(5):2497--2581, 2010.

\bibitem{VP-Sullivan}
Micheline Vigu{\'e}-Poirrier and Dennis Sullivan.
\newblock The homology theory of the closed geodesic problem.
\newblock {\em Journal of Differential Geometry}, 11(4):633--644, 1976.

\bibitem{Viterbo-Functors2}
Claude Viterbo.
\newblock Functors and {C}omputations in {F}loer homology with applications,
  {II}.
\newblock {\em arXiv:1805.01316}.

\bibitem{Viterbo-WeinsteinR2n}
Claude Viterbo.
\newblock A proof of {W}einstein’s conjecture in $\mathbb{R}^{2n}$.
\newblock {\em Annales de l'Institut Henri Poincare (C) Non Linear Analysis},
  4(4):337--356, 1987.

\bibitem{Viterbo-ICM}
Claude Viterbo.
\newblock Generating functions, symplectic geometry, and applications.
\newblock In {\em Proceedings of the International Congress of Mathematicians},
  pages 537--547. Springer, 1995.

\bibitem{Viterbo-Functors1}
Claude Viterbo.
\newblock Functors and computations in {F}loer homology with applications, {I}.
\newblock {\em Geometric \& Functional Analysis}, 9(5):985--1033, 1999.

\bibitem{Weinstein-Conjecture}
Alan Weinstein.
\newblock On the hypotheses of {R}abinowitz' periodic orbit theorems.
\newblock {\em Journal of differential equations}, 33(3):353--358, 1979.

\end{thebibliography}

\end{document}